\documentclass[10pt]{amsart}
\usepackage{bbm}
\usepackage{txfonts}
\usepackage{amsfonts}
\usepackage{mathrsfs}
\usepackage{amssymb}
\usepackage{amsmath}
\usepackage{epic}

\newtheorem{proposition}{Proposition}[section]
\newtheorem{lemma}[proposition]{Lemma}
\newtheorem{corollary}[proposition]{Corollary}
\newtheorem{theorem}[proposition]{Theorem}

\theoremstyle{definition}

\newtheorem{example}[proposition]{Example}

\theoremstyle{remark}
\newtheorem{remark}[proposition]{Remark}

\newcommand{\thlabel}[1]{\label{th:#1}}
\newcommand{\thref}[1]{Theorem~\ref{th:#1}}
\newcommand{\selabel}[1]{\label{se:#1}}
\newcommand{\seref}[1]{Section~\ref{se:#1}}
\newcommand{\lelabel}[1]{\label{le:#1}}
\newcommand{\leref}[1]{Lemma~\ref{le:#1}}
\newcommand{\prlabel}[1]{\label{pr:#1}}
\newcommand{\prref}[1]{Proposition~\ref{pr:#1}}
\newcommand{\colabel}[1]{\label{co:#1}}
\newcommand{\coref}[1]{Corollary~\ref{co:#1}}
\newcommand{\relabel}[1]{\label{re:#1}}

\newcommand{\eqlabel}[1]{\label{eq:#1}}
\newcommand{\equref}[1]{(\ref{eq:#1})}
\def\a{\alpha}
\def\b{\beta}

\def\cd{\cdot}
\def\cds{\cdots}
\def\d{\delta}
\def\D{\Delta}

\def\ep{\varepsilon}
\def\End{\mathrm{End}}
\def\g{\gamma}
\def\G{\Gamma}
\def\Hom{\mathrm{Hom}}

\def\l{\lambda}

\def\lra{\longrightarrow}

\def\mcd{\mathcal{D}}
\def\mcf{\mathcal{F}}

\def\mcm{\mathcal{M}}

\def\mcw{\mathcal{W}}

\def\ml{M(\l)}

\def\N{\mathbb{N}}
\def\op{\oplus}
\def\ot{\otimes}
\def\om{\omega}

\def\ov{\overline}

\def\qgg{\mathcal{U}^{\>0}}
\def\ra{\rightarrow}

\def\su{\mathfrak{u}}
\def\t{\theta}
\def\ti{\times}

\def\ud{U(\mathcal{D},\l)}

\def\vf{\varphi}
\def\vp{\varpi}

\def\wt{\widetilde}
\def\yd{\mathcal {Y}\mathcal {D}}
\def\Z{\mathbb{Z}}

\def\<{\leq}
\def\>{\geq}

\newcommand{\mat}[2]{\begin{bmatrix}
                      #1 \\
                      #2 \\
                    \end{bmatrix}
}

\date{}
\begin{document}
\title{Representation theory of a Class of Hopf algebras}
\author{Zhen Wang and Huixiang Chen}
\address{Department of Mathematics, Yangzhou University,
Yangzhou 225002, China} \email{wangzhen118@gmail.com,
hxchen@yzu.edu.cn} \subjclass{16W30, 17D37} \keywords{pointed Hopf
algebra, quantum group, representation}

\begin{abstract}
The representations of the pointed Hopf algebras $U$ and $\su$ are
described, where $U$ and $\su$ can be regarded as deformations of
the usual quantized enveloping algebras $U_q(\mathfrak{sl}(3))$ and
the small quantum groups respectively. It is illustrated that these
representations have a close connection with those of the quantized
enveloping algebras $U_q(\mathfrak{sl}(2))$ and those of the half
quantum groups of $\mathfrak{sl}(3)$.
\end{abstract}
\maketitle

\section*{\bf Introduction}

In a series of papers \cite{and1,and2, and4,and3}, N.Andruskiewitsch
and H.-J. Schneider classified finite dimensional pointed Hopf
algebras over algebraically closed field $k$ with char$(k)=0$, whose
group-like elements form a finite abelian group. Under some suitable
hypotheses, given a datum of finite Cartan type $\mcd=\mcd
(\G,(g_i)_{1\< i\< \t},(\chi _i)_{1\< i\< \t },(a_{ij})_{1\< i,j\<\t
}),$ a linking datum $\l =(\l _{ij})_{1\<i<j\<\t , i\nsim j}$ and a
root datum $\mu=(\mu_{\a})_{\a \in R^+}$, one can construct a
pointed Hopf algebra $\ud$ and a finite dimensional quotient
$\su(\mcd,\l,\mu)$. On the other hand, if $A$ is a finite
dimensional pointed Hopf algebra over $k$ and the group-like
elements in $A$ form an abelian group, then $A\simeq
\su(\mcd,\l,\mu)$ for some data $\mcd, \l,\mu$. The Hopf algebras
$\ud$ (resp. $\su(\mcd,\l,\mu)$) have close connection with the
quantized enveloping algebras and can be regarded as generalization
of the quantized enveloping algebras (resp. the small quantum
groups).

It is natural to investigate the properties and the representations
of these Hopf algebras. However, not very much is known about this
problem in general. A classical example discussed intensely is the
representations of the quantized enveloping algebras of semisimple
Lie algebras and their quotients (i.e., the small quantum groups).
Recently, the irreducible representations of a class of generalized
doubles are described in \cite{ras}, which can be parameterized by
dominant pairs of characters. Another example is the representation
theory of the half quantum group. In \cite{wch}, simple modules,
projective weight modules and simple Yetter-Drinfeld weight modules
over the half quantum group are fully described.

This paper aims to study the representations of the pointed Hopf
algebras $U=\ud$ and its finite dimensional quotient Hopf algebras
$\su$, where $\mcd $ consists of a free abelian group $\G$ with rank
$2$, elements $(g_1,g_2,g_3,g_4)$ in $\G$, the characters
$(\chi_1,\chi_2,\chi_3,\chi_4)$ and the Cartan matrix of type
$A_2\ti A_2$. In particular, we concentrate on the case that the
linking datum $\l=(\l _{ij})_{1\<i<j\<\t , i\nsim j}$ is given by
$\l_{13}=1$ and $\l_{ij}=0$ otherwise. We describe the simple
modules over $U$ and the simple modules and the projective modules
over $\su$. We note that the category $\mathcal{F}$ of all finite
dimensional $U$-modules is not semisimple. That is different from
the representation theory of the quantized enveloping algebras. In
fact, the representation theory of $U$ can be regarded as
``combination'' of those of the quantized enveloping algebra
$U_q(\mathfrak{sl}(2))$ and those of the half quantum groups of
$\mathfrak{sl}(3)$.

In \seref{1}, some definitions and the structure of the Hopf
algebras $\ud$ and $\su(\mcd,\l,\mu)$ are given. In \seref{2}, Hopf
algebra $U=\ud$ is constructed, where $\mcd$ is a Cartan datum of
type $A_2\ti A_2$. The representation theory of $U$ is also
developed when $q$ is not a root of unity. When $q$ is a root of
unity, there is a Hopf ideal in $U$. The corresponding quotient Hopf
algebra $\su$ is defined in \seref{3}. For a given skew pairing
$\vf$, the double crossproduct $D_{\vf}=D_{\vf}(\su^{\>0},\su
^{\<0})$ is constructed, which is a twisting of the usual Drinfeld
quantum double $D(\su^{\>0},\su ^{\<0})$. In \thref{239}, we show
that the category $_{D_{\vf}} \mcm$ is equivalent to the direct
product of $|\zeta_{l}\ti \zeta_l|$ copies of the category $_{\su}
\mcm$. The simple modules and projective modules over $\su$ are
constructed in Section $3.2$. In particular, we give an equation set
which can be used to compute idempotent elements in $\su$. Thus we
get a decomposition of the regular module as a direct sum of
indecomposable projective modules. This method is valid for the
small quantum group of $\mathfrak{sl}(2)$ too.

Throughout, we assume that $k$ is an algebraically closed field with
characteristic $0$, and all vector spaces and tensor products are
over $k$. Let $k^{\ti}=k\backslash \{0\}$, $\Z$ be the integer set,
$\N$ the positive integer set and $\N_0=\N\cup \{0\}$. For $l\in
\Z$, let $\Z_l=\Z/l\Z$.

\section{\bf Preliminaries}\selabel{1}

For $n\in \N$ and $q\in k\backslash \{0,\pm 1\}$, let
$[n]_q=(q^{n}-q^{-n})/(q-q^{-1})$. As usual, we define $[0]_q!=1$
and $[n]_q!=[n]_q[n-1]_q\cdots \cdots [1]_q$ for $n\geq 1$, and the
Gaussian $q$-binomial coefficients
$$\left [ \begin{array}{c}n\\ j \end{array} \right ]_{q}
=\frac{[n]_q!}{[j]_q![n-j]_q!}\ , \hspace{1cm} n\geq j\geq 0.$$

For $a\in k$, let $[a;n]_q=(aq^{n}-a^{-1}q^{-n})/(q-q^{-1})$.
Similarly we can define the factorial $[a;0]_q!=1,$
$[a;n]_q!=[a;n]_q[a;n-1]_q\cdots \cdots [a;1]_q$ for $n\>1$ and the
``binomial coefficients''
$$\left [ \begin{array}{c}a;n\\ j \end{array} \right
]_{q} =\frac{[a;n]_q!}{[j]_q![a;n-j]_q!}\ , \hspace{1cm} n\geq j\geq
0.$$ Note that $\left [ \begin{array}{c}aq^{-1};n+1\\ j \end{array}
\right ]_{q}=\left [ \begin{array}{c}a;n\\ j \end{array} \right
]_{q}$. For $j\>0$, one can define $$\left [
\begin{array}{c}a;n\\ j \end{array} \right ]_{q}=
\left [ \begin{array}{c}aq^{-j-s};n+j+s\\ j \end{array} \right
]_{q}\hspace{1cm} \text{if }n<j  ,$$ where $s\in \Z$ satisfies
$n+s\>0$.

Let $A, H$ be bialgebras. A bilinear form $\vf :A\ot H\lra k$ is
skew pairing if the following conditions are satisfied:
\begin{eqnarray*}
&\vf(ab,x)=\sum \vf (a,x_{(1)})\vf(b,x_{(2)}),\\
&\vf(a,xy)=\sum \vf (a_{(1)},y)\vf(a_{(2)},x),\\
&\vf(a,1)=\ep(a),\\
&\vf(1,x)=\ep(x),
\end{eqnarray*}
for all $a,b\in A,x,y\in H$. If $A$ (resp. $H^{op}$) is Hopf algebra
with antipode $S_A$ (resp. $S_{H^{op}}$), then $\vf$ is invertible
with $\vf^{-1}(a,x)=\vf(S_A(a),x)\,( resp.
\,=\vf(a,S_{H^{op}}(x)))$. Hence if $A$ and $H$ are Hopf algebras
and $S_H$ is bijective, $\vf(S_A(a),x) =\vf(a,S_{H^{op}}(x))$, where
$S_{H^{op}}=S_H^{-1}$, the composition inverse of $S_H$.

If $\vf$ is a convolution invertible skew pairing, then the double
crossproduct $A\bowtie _{\vf}H=A\ot H$ is constructed as follows.
The coalgebra structure is given by
\begin{eqnarray*}
&&\D (a\ot x)=\sum a_{(1)}\ot x_{(1)}\ot a_{(2)}\ot x_{(2)},\\
&&\ep (a\ot x)=\ep(a)\ep (x),
\end{eqnarray*}
and the algebra structure is given by
\begin{eqnarray*}
&&(a\ot x)(b\ot y)=\sum \vf (b_{(1)}, x_{(1)})ab_{(2)}\ot
x_{(2)}y \vf^{-1}(b_{(3)}, x_{(3)}), \\
&&\text{with identity }1\ot 1.
\end{eqnarray*}
If $A$ and $H$ are Hopf algebras, then $A\bowtie _{\vf}H=A\ot H$ is
also a Hopf algebra (see \cite{doi, maj,rad}).

Recall that a datum of Cartan type $\mcd =\mcd (\G,(g_i)_{1\< i\<
\theta},(\chi _i)_{1\< i\< \theta },(a_{ij})_{1\< i,j\<\theta })$
consists of an abelian group $\G $, elements $g_i\in \G $, $\chi
_i\in \hat{\G }$, $1\< i\<\theta $, and a generalized Cartan matrix
$(a_{ij}) $ of size $\theta $ satisfying the relations:

\begin{equation}\label{1a}
q_{ij}q_{ji}=q_{ii}^{a_{ij}},\, q_{ii}\neq 1, \text{ where }
q_{ij}=\chi _j(g_i),\, 1\< i,j\<\theta .
\end{equation}

$\theta$ is called the rank of $\mcd $. If the matrix $(a_{ij})$ is
of finite type, then $\mcd$ is said to be of finite Cartan type. In
this case, $(a_{ij})$ is a matrix of blocks corresponding to the
connected components of the Dynkin diagram after a reordering. We
write $i\sim j$ for any $i,j\in \{1,\cds,\theta\}$ if $i$ and $j$
are in the same connected component. Let $\mathcal
{I}=\{I_1,I_2,\cds,I_t\}$ be the set of connected components of
$I=\{1,\cds,\t \}$. Let $ord(q)$ denote the order of $q$. For any
$1\< i\<\t$, when $ord(q_{ii})$ is finite, we assume
\begin{eqnarray}
&&ord(q_{ii}) \text{ is odd, and} \label{2}\\
&&ord(q_{ii}) \text{ is prime to }3,\text{ if }i \text{ lies in a
component }G_2.\label{3}
\end{eqnarray}

\begin{lemma}
Given a $J\in \mathcal{I}$, then $ord(q_{ii})=ord(q_{jj})$ for any
$i,j \in J$. We write $N_J=ord(q_{ii})$ for some $i\in J$.
\end{lemma}
\begin{proof}
If each $ord (q_{ii})<\infty$, the claim follows from \cite[Lemma
2.3]{and1}. Now we assume that there is an $i\in J$ such that
$ord(q_{ii})=\infty $. Let $j\in J$. Since the Dynkin diagram of $J$
is connected, there is a chain $i=i_1-i_2-\cds-i_p=j$ in $J$ such
that $a_{i_si_{s+1}}\neq 0$ for $s=1,2,\cds ,p-1$. By (\ref{1a}),
$q_{ii}^{a_{ij}}=q_{jj}^{a_{ji}}$. Thus we have
$ord(q_{jj})=\infty$.
\end{proof}

A family $\l =(\l _{ij})_{1\<i<j\<\t , i\nsim j}$ of elements in $k$
is called a family of linking parameters for $\mcd$ if the following
condition is satisfied for all $1\<i<j\<\t$ with $i\nsim j$,
\begin{eqnarray}\label{4}
\text{if } g_ig_j=1 \text{ or }\chi _i\chi _j\neq \ep, \text{ then }
\l_{ij}=0.
\end{eqnarray}

Vertices $i,j$ are called linkable if $i\nsim j,\,g_ig_j\neq 1 ,\,
\chi _i\chi _j= \ep$. $i,j$ are called linked if $\l_{ij}\neq 0$.
For convenience, let $\l_{ji}=-q_{ji} \l_{ij}$ for all
$1\<i<j\<\t,\,i\nsim j$.

We collect some useful facts from \cite{and2, did}: any vertex $i$
is linkable to at most one vertex $j$; if $i,j$ are linkable, then
$q_{ii}=q_{jj},\, q_{ij}q_{ji}=1$; if $i$ and $k$, respectively, $j$
and $l$, are linkable, then $a_{ij}=a_{kl}, \, a_{ji}=a_{lk}$.

For a given datum $\mcd =\mcd (\G,(g_i)_{1\< i\< \theta},(\chi
_i)_{1\< i\< \theta },(a_{ij})_{1\< i,j\<\theta })$ of finite Cartan
type and a linking datum $\l =(\l _{ij})_{1\<i<j\<\t , i\nsim j}$,
one can construct a Hopf algebra $\ud$ as follows.

Let $V$ be a vector space with a basis $\{x_1,x_2,\cds ,x_{\t}\}$.
Then $V$ or $(V,\cd,\rho)$ is a (left-left) Yetter-Drinfeld module
over $k\G$ such that $x _i \in V_{g_i}^{\chi_i}$, i.e., $g\cd
x_i=\chi_i(g)x_i$ for all $g\in \G$ and $\rho (x_i)=g_i\ot x_i$. The
corresponding braiding $c:V\ot V\ra V\ot V$ is given by $c(x_i\ot
x_j)=q_{ij}x_j\ot x_i$, $1\<i,j\<\t$. Let $_{\G}^{\G}\yd$ denote the
category of all the (left-left) Yetter-Drinfeld modules over $k\G$.
Then $_{\G}^{\G}\yd$ is a braided monoidal category and the free
algebra $k\langle x_1,\cds,x_{\t}\rangle$ is a braided Hopf algebra
in $_{\G}^{\G}\yd$. So the smash product $k\langle
x_1,\cds,x_{\t}\rangle\# k\G $ is a usual Hopf algebra. Let $\ud$ be
the quotient Hopf algebra of $k\langle x_1,\cds,x_{\t}\rangle\# k\G
$ modulo the ideal generated by the elements:
\begin{eqnarray}\label{11a}
(ad_cx_i)^{1-a_{ij}}(x_j),\text{    for all }1\<i,j\<\t,\,i\sim
j,\,i\neq j,\\ \label{11b}
x_ix_j-q_{ij}x_jx_i-\l_{ij}(1-g_ig_j),\text{    for all
}1\<i<j\<\t,\,i\nsim j.
\end{eqnarray}
where $(ad_cx_i)(y)=x_iy-q_{ij_1}q_{ij_2}\cds q_{ij_s}yx_i\eqqcolon
 [x_i,y]_c$ for $y=x_{j_1}x_{j_2}\cds x_{j_s}$, $[x_i,y]_c$ is called
 the braided commutator.

The coalgebra structure of $\ud$ is given by
$$\D (x_i)=g_i\ot x_i+x_i\ot 1, \, \D (g)=g\ot g,\text{ for all }
1\<i\<\t ,g\in \G.$$

For a given indecomposable Cartan matrix $(c_{ij})_{n\ti n}$, let
$P=\sum_{i=1}^n\Z \vp_i$ be the weight lattice. Define simple roots
by
$$ \a_j=\sum_{i=1}^n c_{ij}\vp_i, \quad
j=1,\cdots,n.$$

Let $\D=\{\a_1,\cdots ,\a_n\}$, $Q=\Z\D$ (the root lattice), and
$Q_+=\sum_i \N_0\a_i$. For any $\b=\sum_i b_i\a_i\in Q$, define
$g_{\b}=g_1^{b_1}g_2^{b_2}\cds g_n^{b_n}$. In particular,
$g_{\a_i}=g_i$.

Define automorphisms $\g _i$ of $P$ by $\g _i\vp _j=\vp _j-\d_{ij}\a
_i\, (i,j=1,\cdots,n).$ Then $\g_i \a_j=\a_j-c_{ij}\a_i$. Let $W$ be
the (finite) subgroup of $GL(P)$ generated by $\g_1,\cdots, \g_n$,
called the Weyl group. Then $Q$ is $W$-invariant. Let
$R=W\D,\,R^+=R\cap Q_+$ and $R^-=-R^+$. Then $R$ is a root system
corresponding to the Cartan matrix $(c_{ij})$, $R^+$ the set of
positive roots, $R=R^+\cup R^-$.

Fix a reduced expression $\g_{i_1}\g_{i_2}\cdots \g_{i_N}$ of the
longest element $\om _0$ of $W$. This gives us a convex ordering of
the set of positive roots $R^+$:
$$\b _1=\a_{i_1},\, \b_2=\g_{i_1}\a_{i_2},\,\cdots,\,
\b_N=\g_{i_1}\cdots \g_{i_{N-1}}\a_{i_N}.$$

Then for any $J\in \mathcal{I}$ one can choose a Weyl group $W_J$
and a root system $R_J$ for $(a_{ij})_{i,j\in J}$ and a reduced
expression of the longest element $w_{0,J}$ in the Weyl group $W_J$.
Put
$$w _0\coloneqq
w_{0,I_1}w_{0,I_2}\cds w_{0,I_t}$$ and
$$R
^+\coloneqq \{\b_{I_1,1},\cds,\b_{I_1,p_{I_1}},\cds,
\b_{I_t,1},\cds,\b_{I_t,p_{I_t}}\},$$ where $p_J$ is the number of
the positive roots in $R_J^+$ and $\b_{J,1},\cds,\b_{J,p_J} \in
R_J^+$ with the convex ordering. We also write $$R
^+=\{\b_1,\cds,\b_p\},\quad p=\sum_{J\in \mathcal{I}}p_J$$ with the
given ordering.

For each $\b_i \in R_J^+ $, one can define a root vector
$x_{\b_i}\in \ud$ by the same iterated braided commutator of the
elements $x_j,j\in J$ as the Lusztig's case in \cite{lus} but with
respect to the general braiding $c$, see \cite{and4}.

The following theorem describes the structure of $\ud$, which was
stated in \cite[Theorem~3.3]{and1} for a finite abelian group $\G$.
Indeed, the finiteness condition of $\G$ is not necessary in the
proof of \cite[Theorem~3.3]{and1}.
\begin{theorem}\thlabel{111}
Let $\G$ be an abelian group, and $\mcd$ a datum of finite Cartan
type satisfying the conditions (\ref{2}) and (\ref{3}). Let $\l$ be
a family of linking parameters for $\mcd$. Then
\begin{enumerate}
              \item The elements$$x_{\b_1}^{a_1}x_{\b_2}^{a_2}\cds
              x_{\b_p}^{a_p}g, \, a_1,a_2,\cds, a_p\>0,g\in \G,$$
              form a basis of the vector space $\ud$.
              \item Let $J\in \mathcal{I}$ and $\a \in R ^+,\b
              \in R _J^+$. Then $[x_{\a},x_{\b}^{N_J}]_c=0$, that
              is
              $$x_{\a}x_{\b}^{N_J}=q_{\a\b}^{N_J}x_{\b}^{N_J}x_{\a}.$$
            \end{enumerate}

\end{theorem}

A family $\mu=(\mu_{\a})_{\a \in R^+}$ of elements in $k$ is called
a family of root vector parameters for $\mcd $ if the following
condition is satisfied for all $\a \in R_J^+$, $J\in \mathcal{I}$:
If $g_{\a}^{N_J}=1$ or $\chi _{\a}^{N_J}\neq \ep$, then
$\mu_{\a}=0$.

Then we define $$u(\mcd ,\l, \mu )=U(\mcd,
\l)/(x_{\a}^{N_J}-u_{\a}(\mu)|\a \in R_J^+, J\in \mathcal{I}), $$
where $u_{\a}(\mu)$ is central in $U(\mcd, \l)$ and is determined by
$\mu$ uniquely (see \cite{and1}).

\begin{theorem}\thlabel{1.20}(\cite[Thm6.2]{and1})
Let $A$ be a finite dimensional pointed Hopf algebra with abelian
group $G(A)=\G$ and infinitesimal braiding matrix $(q_{ij})_{1\<i,
j\<\t}$. Assume that the following conditions are satisfied:
\begin{align*}
    &ord(q_{ii})>7\text{ is odd},\\
    &ord(q_{ii})\text{ is prime to }3\text{ if } q_{il}q_{li}\in
    \{q_{ii}^{-3},q_{ll}^{-3}\}\text{ for some }l,
\end{align*}
where $1\<i\<\t$. Then
$$A\cong u(\mcd ,\l, \mu ),$$
where $\mcd =\mcd (\G, (g_i)_{1\<i\<\t},
(\chi_i)_{1\<i\<\t},(a_{ij})_{1\<i,j\<\t})$ is a datum of finite
Cartan type, and $\l$ and $\mu $ are families of linking and root
vector parameters for $\mcd$. \end{theorem}

We apply these construction to the special case that the Cartan
matrix $(a_{ij})$ is of type $A_2 \times A_2$, i.e.
$$(a_{ij})=\left(
                                             \begin{array}{cccc}
                                               2 & -1 & 0 & 0 \\
                                               -1 & 2 & 0 & 0 \\
                                               0 & 0 & 2 & -1 \\
                                               0 & 0 & -1 & 2 \\
                                             \end{array}
                                           \right).$$
Let $\G =\langle g_1,g_2\rangle $ be a free abelian group of rank
$2$, $g_3=g_1,\,g_4=g_2$. Let $q\in k\backslash \{0,\pm 1\}$. Then
$\chi _j$ is given by $\chi_j (g_i)=q^{a_{ij}}$ for $1\<j,i\<2$, and
$\chi _3=\chi_1^{-1},\, \chi_4=\chi _2^{-1}$. These form a datum of
finite Cartan type $\mcd$.

The linking datum $(\l_{ij})$ has the following $3$ cases after a
suitable permutation:
$$\left(
 \begin{array}{cccc}
  0 & 0 & 0 & 0 \\
  0 & 0 & 0 & 0 \\
   0 & 0 & 0 & 0 \\
   0 & 0 & 0 & 0 \\
   \end{array}
   \right),\left(
 \begin{array}{cccc}
  0& 0 & 1 & 0 \\
  0 & 0 & 0 & 0 \\
   -q_{31} & 0 & 0 & 0 \\
   0 & 0 & 0 &0 \\
   \end{array}
   \right),
   \left(
 \begin{array}{cccc}
 0 & 0 & 1 & 0 \\
  0 & 0 & 0 & 1\\
   -q_{31} & 0 &0 & 0 \\
   0 & -q_{42} &0 & 0 \\
   \end{array}
   \right),    $$
where $q_{31}=q^{-a_{31}},q_{42}=q^{-a_{42}}$. Denote them by
$\l_1,\l_2,\l_3$ respectively. Then we have the following facts:
 \begin{itemize}
 \item $U(\mcd ,\l_1)$ is a graded Hopf algebra if we define
 $\text{deg}g=0,\,\forall g\in \G ,\,degx_i=1,\, 1\<i\<4$. This is
 similar to the half quantum group discussed in \cite{cib, wch}.\\
  \item $U(\mcd ,\l_3)$ is the quantized enveloping algebra of
 $\mathfrak{sl}(3)$, which has been discussed in many papers, see \cite{brg,dck,jan}.
 \end{itemize}
In the next section, we shall concentrate on the second case, i.e.,
$U(\mcd ,\l_2)$.

\section{\bf Hopf algebra $U(\mcd ,\l_2)$}\selabel{2}
\subsection{the properties of $U$}\selabel{2.1}

Let $\G =\langle g_1,g_2\rangle $ be a free abelian group of rank
$2$, $g_3=g_1,\,g_4=g_2$. Let $q\in k\backslash \{0,\pm 1\}$ and
$(a_{ij})$ the Cartan matrix of type $A_2$. Then $\chi _j$ is given
by $\chi_j (g_i)=q^{a_{ij}}$ for $1\<j,i\<2$. Let $\chi
_3=\chi_1^{-1}$ and $ \chi_4=\chi _2^{-1}$. These form a datum of
finite Cartan type $\mcd$, the corresponding Cartan matrix is of
type $A_2\times A_2$.  For any $1\< i\<4$, if $ord(q_{ii})$ is
finite, we assume that the condition (2) is satisfied.

For convenience, write $U'$ for the corresponding Hopf algebra
$U(\mcd ,\l_2)$.

As an algebra, $U'$ is generated by the elements $x_i(1\<i\<4)$,
$g\in \G$ and subjects to the relations $(\text{for all
}1\<i,j\<4)$:
\begin{eqnarray}
&&gx_j=\chi _j(g)x_jg ,\\
&&ad_c(x_i)^{1-a_{ij}}(x_j)=0,\text{    for all }i\sim
j,\,i\neq j,\\
&&x_ix_j-q_{ij}x_jx_i=\l_{ij}(1-g_ig_j),\text{    for all }i\nsim j.
\end{eqnarray}
where $ad_c(x_i)(x_j)=x_ix_j-q^{a_{ij}}x_jx_i$. In the following we
define an algebra $U$ which is similar to the usual quantized
enveloping algebra $U_q(\mathfrak{sl}(3))$.

Let $U$ be an algebra generated by $E_i,\,F_i$ and $K_i^{\pm
1}(1\<i\<2)$ satisfying the relations:

\begin{eqnarray}
&&K_iK_j=K_jK_i,\, K_iK_i^{-1}=K_i^{-1}K_i=1, \\
&&K_iE_jK_i^{-1}=q^{a_{ji}}E_j, \\
&&K_iF_jK_i^{-1}=q^{-a_{ji}}F_j,\\
&&E_iF_j-F_jE_i=\d_{ij}\d_{i1}(K_i-K_i^{-1})/(q-q^{-1}),\\
&&\sum_{s=0}^{2}(-1)^s\left [ \begin{array}{c}2\\ s
\end{array} \right ]_qE_i^{2-s}E_jE_i^s=0 \quad \text{if }i\neq j,\\
&&\sum_{s=0}^{2}(-1)^s\left [ \begin{array}{c}2\\ s
\end{array} \right ]_qF_i^{2-s}F_jF_i^s=0 \quad \text{if }i\neq j.
\end{eqnarray}

The comultiplication $\D$, antipode $S$ and counit $\ep$ of $U $ are
defined by$(i=1,2)$
\begin{eqnarray*}
&\D(E_i)=K_i\ot E_i+E_i\ot 1,&\D(F_i)=1\ot F_i+F_i\ot K_i^{-1},\\
&\D(K_i)=K_i\ot K_i,&\D(K_i^{-1})=K_i^{-1}\ot K_i^{-1},\\
&S(E_i)=-K_i^{-1}E_i,&S(F_i)=-F_iK_i,\quad S(K_i)=K_i^{-1},\\
& \ep (E_i)=0=\ep (F_i), &\ep (K_i)=1 .
\end{eqnarray*}

\begin{proposition}\prlabel{}
There is a Hopf algebra isomorphism $\phi: U\lra U'$ such that
$$E_i\mapsto x_i,\, F_i\mapsto
(q^{-1}-q)^{-1}x_{i+2}g_i^{-1},\,K_i\mapsto g_i$$ for all $1\<i\<2$.
\end{proposition}
\begin{proof}
It is straightforward.
\end{proof}

Hence we may discuss the properties and the representation of $U$
instead of $U'$.

Similarly to the discussion of the quantized enveloping algebra, let
$U ^+, \,U ^- ,\,U ^0$ be the subalgebras of $U$ generated by the
$E_i$, the $F_i$, and the $K_i,K_i^{-1}(1\< i\< 2)$ respectively. It
follows from $(10)-(15)$ that $U =U ^pU ^qU ^r$, where $(p,q,r)$ is
a permutation of $(+,-,0)$. Moreover, the multiplication gives a
$k$-vector space isomorphism
\begin{equation}\label{16}
U ^p\ot U ^q\ot U ^r \simeq U,\quad \text{where }(p,q,r) \text{ is
as above}.
\end{equation}
Let $U^{\<0}=U^-U^0$ and $U^{\>0}= U^0U^+$ be the Borel subalgebras
of $U$.

As in Section $1$, we can define the corresponding weight lattice
and root lattice for the Cartan matrix of type $A_2$. Let
$P=\sum_{i=1}^2\Z \vp_i$ be the weight lattice. Define simple roots
as follows:
$$ \a_j=\sum_{i=1}^2 a_{ij}\vp_i, \quad
j=1,2.$$

Let $\D=\{\a_1 ,\a_2\}$ be the set of simple roots and $Q=\Z\D$ be
the root lattice. Then there is a $\Z$-bilinear map $(,)$ on
$P\times Q$ given by $(\vp_i, \a_j)=a_{ij}$. Clearly, $(,)$ is
non-degenerate. Let $Q_+=\sum_i \N_0\a_i$. Then there is a partial
ordering on $P$ defined by $\l \>\mu$ if $\l -\mu \in Q_+$.

Put $E_{1,2}=E_1E_2-q^{-1}E_2E_1$ and $ F_{1,2}=F_2F_1-qF_1F_2$.
Then we have the following relations from $(10)-(15)$:
\begin{eqnarray}
&&K_iE_{1,2}K_i^{-1}=qE_{1,2}, \quad K_iF_{1,2}K_i^{-1}=q^{-1}F_{1,2},\\
&&E_{1,2}E_1=q^{-1}E_1E_{1,2},\quad E_{1,2}E_2=qE_2E_{1,2},\\
&&F_{1,2}F_1=q^{-1}F_1F_{1,2},\quad F_{1,2}F_2=qF_2F_{1,2},\\
&&E_{1,2}F_{1,2}-F_{1,2}E_{1,2}=0,\\
&&E_{1,2}F_1-F_1E_{1,2}=-E_2K_1^{-1},\quad
E_{1,2}F_2-F_2E_{1,2}=0,\\
&&F_{1,2}E_1-E_1F_{1,2}=K_1F_2,\quad F_{1,2}E_2-E_2F_{1,2}=0.
\end{eqnarray}

For $\textbf{s}=(s_1,s_2,s_3)\in \N_0^3,\,\l =\sum_{i=1}^2 r_i\a_i
\in Q$, let $E^{\textbf{s}}=E_1^{s_1} E_{12}^{s_2}E_2^{s_3}$,
$F^{\textbf{s}}=F_1^{s_1} F_{12}^{s_2}F_2^{s_3}$ and
$K_{\l}=K_1^{r_1} K_2^{r_2}$. Then the PBW basis of $U^+ (resp.
U^-)$ is given by $\{ E^{\textbf{s}}:\textbf{s}\in \N_0^3\} (resp.
\{ F^{\textbf{t}}:\textbf{t}\in \N_0^3\})$.

Then $U$ has a basis
$\{F_1^{t_1}F_{1,2}^{t_2}F_2^{t_3}K_1^{r_1}K_2^{r_2}E_1^{s_1}E_{1,2}^{s_2}E_2^{s_3}
|s_i,t_i \in \N_0,r_j \in \Z,1\<i\<3,1\<j\<2\}$.

\subsection{Representations of $U(\mcd ,\l_2)$}

Let $\hat{\G}$ denote the character group of $\G$. Then
$\hat{\G}\simeq (k^{\ti})^2$ as groups. We always identity the
characters with the elements in $(k^{\ti})^2$ below. For a weight
$\l=m_1\vp _1+m_2\vp_2\in P$, define a character $i(\l) \in
\hat{\G}$ by $i(\l)(K_i)=q^{(\l, \a_i)}$. Then $i:P\lra \hat {\G}$
is a group homomorphism and the image of $i$ is a subgroup of
$\hat{\G}$. Moreover, $i$ is a monomorphism if and only if $q$ is
not a root of unity. Define $P_1=\{\l=m_1\vp _1+m_2\vp_2\in
P|0\<m_1,m_2\<l-1\}$ when $q$ is an $l$-th primitive root of unity
and $P_1=P$ otherwise. Then $i(P_1)=i(P)$ and $i|_{P_1}$ is an
injective map. In the sequel, we usually regard $\l \in P_1\subseteq
P$ for some character $\l\in \hat {\G}$ if $\l \in i(P)$.

In the rest of this section, we assume that $q$ is not a root of
unity.

For any character $\l \in \hat{\G}$, let $M(\l)\coloneqq U\ot
_{U^{\>0}}V$, where $V=kv$ is a $U^{\>0}$-module with the action
given by $E_i\cd v =0$ and $K_i\cd v=\l(K_i)v$, $i=1,2$. Then
$M(\l)= \text{span}\{F_1^{t_1}F_{1,2}^{t_2}F_2^{t_3}\cd v|t_i \in
\N_0\}$ as vector spaces, where we identity $1\ot v$ with $v$. That
is, $M(\l)$ is a free module of rank $1$ over $U^-$.

Note that $E_2\cd F_1^{t_1}F_{1,2}^{t_2}F_2^{t_3}\cd v=0$. Let
$$M[\l, n]\coloneqq \text{span}\{F_1^{t_1}F_{1,2}^{t_2}F_2^{t_3}\cd v| t_i\>0, t_2+t_3\>n\}
\quad \text{for all } n\>0.$$

\begin{proposition}\prlabel{}
For any $n\>0$, $M[\l,n]$ is a submodule of $M(\l)$.
\end{proposition}
\begin{proof}
It is a straightforward verification.
\end{proof}
Then there is a filtration
\begin{equation}\eqlabel{2a}
 M(\l)=M[\l,0]\supseteq
M[\l, 1]\supseteq M[\l,2]\supseteq \cds .
\end{equation}

Let $M(\l,n)=M[\l,n]/M[\l,n+1]$. Let $\pi _n:M(\l)\lra M(\l)/ M[\l,
n+1]$ be the natural projection. Then $M(\l,n)=\text{span}
\{F_1^{t_1}F_{1,2}^{t_2}F_2^{t_3}\cd \pi_n(v)| t_i\>0, t_2+t_3=n\}$.

\begin{proposition}\prlabel{}
$F_2 \cd M(\l, n)=0 $ for all $n\>0$.
\end{proposition}
\begin{proof}
It is obvious.
\end{proof}
\begin{remark}
Let $\rho _n :U \ra \End(M(\l,n))$ be the corresponding
representation of $U$. Then $E_2, F_2\in \text{Ker} ~\rho _n$. There
is a natural induced representation $\bar{\rho}_n$ of the quotient
algebra $U/\langle E_2,F_2\rangle$ on $M(\l,n)$.
\end{remark}

Let $M$ be a $U$-module and $\l \in \hat{\G}$. $0\neq v\in M$ is
called weight vector with weight $\l$ if $g\cd v=\l(g)v$ for all
$g\in \G$. Write $M_{\l}$ for the set of all the weight vectors in
$M$ with weight $\l$. Let $\Pi (M)$ be the set of all weights $\l$
with $M_{\l}\neq 0$. Call $M$ a weight module if $M=\op_{\l}M_{\l}$.
A weight vector $ v\in M$ is called highest weight vector if
$E_1v=E_2v=0$. A $U$-module is a highest weight module of weight
$\l$ if it is generated by a highest weight vector with weight $\l$.
Let $\mcw$ denote the category of all weight modules. Clearly, every
highest weight module lies in $\mcw$.

For $n\>0$, let $$v_n\coloneqq a_0F_{1,2}^n\cd
v+a_1F_1F_2F_{1,2}^{n-1}\cd v+\cds +a_nF_1^nF_2^n \cd v\in M(\l),$$
where $a_i=q^{-(n-i)}\l(K_1)^{-(n-i)}[n-i]_q!\begin{bmatrix}
                       n \\
                       n-i\\
                     \end{bmatrix}_q\begin{bmatrix}
                       \l(K_1);1 \\
                       n-i\\
                     \end{bmatrix}_q$ for $0\<i\<n$.
Then $v_0=v$.
\begin{proposition}
For all $n\>0$, $v_n$ defined as above are highest weight vectors.
\end{proposition}
\begin{proof}
This is similar to the proof of \cite[Theorem~1]{wan}.
\end{proof}

Then for any two fixed integers $t_2,t_3\in \N_0$ with $t_2+t_3=n$,
let $M(\l,n;t_2,t_3)=\text{span} \{F_1^{t_1}F_2^{t_3}\cd
\pi_n(v_{t_2})| t_1\>0\}\subseteq M(\l,n)$. It is easy to see that
$M(\l,n;t_2,t_3)$ is a submodule of $M(\l,n)$ and
\begin{align}\label{24}
 M(\l,n)=\bigoplus _{t_2+t_3=n}M(\l,n;t_2,t_3).
\end{align}

Hence $\{F_1^{t_1}F_2^{t_3}v_{t_2}| t_1,t_2,t_3\> 0\}$ is a
$k$-basis of $\ml$.

Let $V$ be a $k$-vector space with a basis $\{u_0,u_1,\cds\}$ and
$\l=(\l_1,\l_2)\in (k^{\ti})^2$. Then one can check that $V$ admits
a $U$-module structure with the $U$-action given by
\begin{align*}
 K_1u_j&=\l_1q^{-2j}u_j,\quad K_2u_j=\l_2q^ju_j,\\
 F_1u_j&=u_{j+1},\quad
 E_1u_j=[j]_q[\l_1;1-j]_qu_{j-1},\\
 E_2u_j&=F_2u_j=0,&
 \end{align*}
where $j=0,1,2,\cds$ and $u_{-1}=0$. Denote the module by $V(\l)$.

\begin{proposition}\prlabel{2.1.5a}
For any $\l ,\mu \in (k^{\ti})^2$, $V(\l)\simeq V(\mu)$ if and only
if $\l =\mu$.
\end{proposition}
\begin{proof}
Let $\{u_i\}(\text{ resp. }\{w_i\})$ be the $k$-basis of
$V(\l)(\text{ resp. }V(\mu))$ on which $U$ acts as above. Let $f:
V(\l)\ra  V(\mu)$ be a non-zero $U$-module homomorphism. Then
$f(u_0)=\sum_{i=0}^r a_iw_i\neq 0$ for some $a_i\in k$,
$i=1,2,\cds,r$. Hence $f(u_j)=\sum_{i=0}^r a_iw_{i+j}\neq 0$. Note
that $f$ is an isomorphism if and only if $\{f(u_j)\}$ is a basis of
$V(\mu)$, if and only if there are some scalars $b_j\in k$,
$j=0,1,\cds, s$, such that $w_0=\sum _{j=0}^s b_jf(u_j)=\sum
_{i,j}a_ib_jw_{i+j}$, if and only if $a_0b_0=1$ and $a_i=b_j=0$ for
$i>0,j>0$, i.e., $f(u_0)=a_0w_0(a_0\neq 0)$. This finishes the proof
by virtue of the actions of $K_1,K_2$.
\end{proof}

\begin{proposition}\prlabel{2.1.6}
Let $\l \in \hat{\G}$ and $M(\l,n;t_2,t_3)$ be defined as above.
Then
$$M(\l,n;t_2,t_3)\simeq V(\l'),$$ where
$\l'=(q^{t_3-t_2}\l(K_1),q^{-t_2-2t_3}\l(K_2))\in (k^{\ti})^2$.
\end{proposition}
\begin{proof}
Define a $k$-map $f: V(\l')\ra M(\l,n;t_2,t_3)$ by
$f(u_i)=F_1^iF_2^{t_3}\cd \pi_n(v_{t_2})$ for all $i\>0$. It is easy
to check that $f$ is a $U$-module isomorphism.
\end{proof}

When $n=0$, then $t_2=t_3=0$. Hence $M(\l,0)=M(\l,0;0,0)\simeq
V(\l)$, where we identity $\l \in \hat{\G}$ with
$\l=(\l(K_1),\l(K_2))\in(k^{\ti})^2$.

Note that these $V(\l)$ are similar to the Verma modules over the
quantized enveloping algebra $U_q(\mathfrak{sl}_2)$.

\begin{proposition}\prlabel{2.1.7}
For $\l=(\l_1,\l_2)\in (k^{\ti})^2$, $V(\l)$ is an indecomposable
$U$-module. Furthermore, $V(\l)$ is reducible if and only if there
exists $m_1\in \N_0$ such that $\l_1=\pm q^{m_1}$.
\end{proposition}
\begin{proof}
Let $V(\l)=V'\op V''$ as $U$-modules. Then $u_0=u'+u''$, where
$u'\in V'$ and $u''\in V''$. We have that
$K_iu'+K_iu''=K_iu_0=\l_iu'+\l_iu''$, $i=1,2$. Since $V',V''$ are
both $U$-modules and $V'\cap V''=0$ , $u'$ and $u''$ are weight
vectors with weight $\l$. Note that $q$ is not a root of unity. Then
$V(\l)_{\l}=ku_0$. This implies that $u'=0$ or $u''=0$, so
$V(\l)=V''$ or $V(\l)=V'$.

If $\l_1=\pm q^{m_1}$ for some $m_1\in \N_0$, then
$J'(\l)=span\{u_{m_1+1},u_{m_1+2},\cds\}$ is a submodule of $V(\l)$.
If not, we claim that $V(\l)$ is irreducible. For any $0\neq
w=\sum_{i=0}^n a_iu_i\in V(\l)$ with $a_n\neq 0$, we have that
$E_1^nw=a_n[n]_q![\l_1;1-n]_q[\l_1;2-n]_q\cds [\l_1;0]_qu_0$. Since
$\l_1\neq \pm q^{m_1}$ for any $m_1\in \N_0$ and $q$ is not a root
of unity, $E_1^nw=au_0$ for a non-zero scalar $a$. Then $Uw=V(\l)$
and the claim follows.
\end{proof}

Obviously, we have

\begin{corollary}
Let $\l \in \hat{\G}$ and $M(\l,n;t_2,t_3)$ be defined as above.
Then $M(\l,n;t_2,t_3)$ is indecomposable. Furthermore,
$M(\l,n;t_2,t_3)$ is reducible if and only if $\l(K_1)=\pm q^{m_1}$
for some $m_1\in \Z$ with $m_1\>t_2-t_3$.
\end{corollary}

Let $\l=(\l_1,\l_2)\in (k^{\ti})^2$. If there is an $m_1\in \N_0$
such that $\l_1=\pm q^{m_1}$, then
$J'(\l)=span\{u_{m_1+1},u_{m_1+2},\cds\}$ is the unique maximal
submodule of $V(\l)$. The induced quotient is an
$(m_1+1)$-dimensional simple module, denoted by $L(\l)$. For
convenience, let $u_i$ denote the image of $u_i$ in $L(\l)$ under
the natural projection, $0\<i\<m_1$. Call the set
$\{u_i:0\<i\<m_1\}$ the standard basis of $L(\l)$.

\begin{proposition}
The finite dimensional simple modules $L(\l)$ are non-isomorphic
each other.
\end{proposition}
\begin{proof}
For two simple modules $L(\l)$ and $L(\mu)$, let $\{u_i:0\<i\<m_1\}$
(resp. $\{w_i:0\<i\<m_1'\}$) be the standard basis of $L(\l)$ (resp.
$L(\mu)$). Assume that there is a $U$-module isomorphism $f:
L(\l)\ra L(\mu)$. Then $\dim L(\l)=\dim L(\mu)=m_1$. Note that
$L(\mu)_{\tau}$ is $1$-dimensional and is spanned by some $w_i$ for
any $\tau \in \Pi (L(\mu))$. Then $f(u_0)=aw_i$ for some $a\in
k^\ti$. From $0\neq f(u_{m_1})=F_1^{m_1}f(u_0)=aF_1^{m_1}w_i$, we
have that $i=0$. Hence $\l=\mu $ by virtue of the actions of
$K_1,K_2$ on $u_0$ and $w_0$.
\end{proof}

In fact, one can show that any finite dimensional simple module over
$U$ must be isomorphic to some $L(\l)$.

Let $J(\l)=M[\l,1]\op \text{span}\{F_1^{t_1}\cd v| t_1\>m_1+1\}$ if
$\l_1=\pm q^{m_1}$ for some $m_1\in \N_0$, and $J(\l)=M[\l,1]$
otherwise.
\begin{proposition}\prlabel{212}
For any $\l \in \hat{\G}=(k^\ti)^2$, $J(\l)$ is the unique maximal
submodule of $M(\l)$.
\end{proposition}
\begin{proof}
It suffices to show that $U\cd w=M(\l)$ for any vector $0\neq w \in
M(\l)\setminus J(\l)$. Since $M(\l)$ is a weight module, one can
assume that $w$ is a weight vector. Write $w=\sum
a_{t_1,t_2,t_3}F_1^{t_1}F_{1,2}^{t_2}F_2^{t_3}\cd v \in M(\l)$ for
$a_{t_1,t_2,t_3}\in k$. Since $q$ is not a root of unity,
$F_1^{t_1}F_{1,2}^{t_2}F_2^{t_3}\cd v$ and
$F_1^{s_1}F_{1,2}^{s_2}F_2^{s_3}\cd v$ have the same weights if and
only if $t_1+t_2=s_1+s_2$ and $t_2+t_3=s_2+s_3$. Then
$w=aF_1^{t_1}\cd v$ for some $a\in k^\ti$ and $t_1\>0$. Furthermore,
$t_1\<m_1$ if $\l_1=\pm q^{m_1}$ for some $m_1\>0$. Then the result
follows from that $E_1^{t_1}F_1^{t_1}\cd
v=[t_1]_q![\l_1;1-t_1]_q[\l_1;2-t_1]_q\cds [\l_1;0]_qv$ for any
$t_1\>1$.
\end{proof}

\begin{proposition}
$M$ is a highest weight module if and only if $M$ is a quotient of
some Verma module $\ml$.
\end{proposition}
\begin{proof}
The part ``if'' is obvious. We need to check the part ``only if''.
Let $M$ be a highest weight module generated by a highest weight
vector $v$ of weight $\l$. Note that $U=U^-U^0U^+=U^-U^0+UE_1+UE_2$.
Then we have that $M=Uv=(U^-U^0+UE_1+UE_2)v=U^-v$ since $v$ is a
highest weight vector. The claim follows from that $\ml$ is a free
$U^-$-module of rank $1$.
\end{proof}

Let $M$ be a highest weight module. Then there exists some
$U$-module $M(\l)$ and some submodule $N\subseteq M(\l)$ such that
$M\cong M(\l)/N$. Regard $M=M(\l)/N$. Recall that $\ml$ has a
filtration \equref{2a}. Let $M[i]\coloneqq  (M[\l,i]+N)/N$ for
$i\>0$. There is an induced filtration of $M$:
\begin{equation*}
M=M[0]\supseteq M[1] \supseteq \cds.
\end{equation*}
Then $M[0]/M[1]$ is isomorphic to some quotient of $V(\l)\simeq
M(\l,0)= M[\l,0]/M[\l,1]$. Hence we have the following result.
\begin{proposition}\prlabel{2112}
Let $M$ be a highest weight module over $U$ with highest weight
$\l$. If $M$ is simple, then
\begin{enumerate}
   \item $M\simeq L(\l)$ if $M$ is finite dimensional;
   \item $M\simeq V(\l)$ if $M$ is infinite dimensional.
 \end{enumerate}
\end{proposition}
\begin{proof}
From the discussion above, one can get a filtration of submodules
$M=M[0]\supseteq M[1] \supseteq \cds.$ Since $M$ is simple, then
$M[1]=0$ or $M[1]=M$. If $M[1]=M$, then $M[\l,1]+N=M(\l)$, where $N$
is given as before. By \prref{212}, $M[\l,1]\subseteq J(\l)$, which
is a small submodule. This implies that $N=M(\l)$ by the Nakayama
Lemma, a contradiction. Thus $M[1]=0$ and $M=M[0]/M[1]$ is
isomorphic to some quotient of $V(\l)$.

If $M$ is finite dimensional, then $\l_1=\pm q^{m_1}$ for some
$m_1\in \N_0$ by \prref{2.1.7}. Then $M\simeq L(\l)$ since $J'(\l)$
is the unique maximal submodule of $V(\l)$ and $L(\l)=V(\l)/J'(\l)$.
If $M$ is infinite dimensional, then $\l_1\neq \pm q^{m_1}$ for any
$m_1\in \N_0$.  By \prref{2.1.7}, $V(\l)$ is irreducible, so
$M\simeq V(\l)$.
\end{proof}

\begin{proposition}\prlabel{213}
Any non-zero finite dimensional $U$-module $M$ contains a highest
weight vector. Moreover, the endomorphisms induced by
$E_1,E_2,F_1,F_2$ are nilpotent.
\end{proposition}
\begin{proof}
Let $M$ be a finite dimensional $U$-module. Since $U^0$ is
commutative and $k$ is algebraically closed, there exists a weight
vector $v \in M$. Let $N\coloneqq Uv$. Then $N$ is a weight module,
write $N=\op _{\l \in(k^{\ti})^2}N_{\l}=\op _{\l \in\Pi(N)}N_{\l}$.
Let $\epsilon_1=(q^2,q^{-1}),\epsilon_2=(q^{-1},q^2)\in k^2$. For a
fixed $\l\in \Pi(N)$, then $E_1\cd N_{\l}\subseteq
N_{\epsilon_1\l},E_2\cd N_{\l}\subseteq N_{\epsilon_2\l}$. Let
$W_0=\{\l\}\subseteq \Pi(N)$ and $\l^{(0)}=\l$. Since $q$ is not a
root of unity and $\Pi(\l)$ is a finite set, there exists an $n_1\in
\N_0$ such that $E_1^{n_1}\cd N_{\l}\neq 0$ and $E_1^{n_1+1}\cd
N_{\l}=0$. Define $W_1=W_0\cup \{\epsilon_1\l,\epsilon_1^2\l,\cds,
\epsilon_1^{n_1}\l\}\subseteq \Pi(N)$ and
$\l^{(1)}=\epsilon_1^{n_1}\l$. Since $E_1^{n_1}\cd N_{\l}\neq 0$,
there exists an $n_2\in \N_0$ such that $E_2^{n_2}\cd
N_{\l^{(1)}}\neq 0$ and $E_2^{n_2+1}\cd N_{\l^{(1)}}=0$. Define
$W_2=W_1\cup \{\epsilon_2\l^{(1)},\epsilon_2^2\l^{(1)},\cds,
\epsilon_2^{n_2}\l^{(1)}\}\subseteq \Pi(N)$ and
$\l^{(2)}=\epsilon_2^{n_2}\l^{(1)}$. Iterating this process, one can
define the sets $W_3,W_4,\cds$, and the weights
$\l^{(3)},\l^{(4)},\cds$. Thus one gets an ascending chain
$W_0\subseteq W_1\subseteq W_2\subseteq \cds $ of subsets of
$\Pi(N)$. Since $q$ is not a root of unity, $W_j=W_{j+1}$ if and
only if $n_{j+1}=0$, where $j\>0$. Moreover, each $N_{\l^{(j)}}\neq
0$.

Since $\Pi(N)$ is a finite set, there exists an $s\in \N_0$ such
that $W_s=W_{s+1}=W_{s+2}$. Without loss of generality, assume that
$s$ is the minimal integer with this property. Since $q$ is not a
root of unity, $W_s=W_{s+1}=W_{s+2}$ implies that
$n_{s+1}=n_{s+2}=0$, that is $E_1\cd N_{\l^{(s)}}=0$ and $E_2\cd
N_{\l^{(s)}}=0$. This shows that there exists a highest weight
vector in $N\subseteq M$.

If $M$ is a simple $U$-module, then $M$ is a highest weight module.
The last claim is obvious. If not, one can take a composition
sequence of $M$. Then the last claim follows from that the
endomorphisms induced by $E_1,E_2,F_1,F_2$ are nilpotent on each
simple factor.
\end{proof}

\begin{remark}
Here the assumption that $k$ is an algebraically closed field is
necessary.
\end{remark}

\begin{corollary}\colabel{213a}
Every finite dimensional simple $U$-module is a highest weight
module.
\end{corollary}

By \prref{2112}, \coref{213a} means that any finite dimensional
simple module over $U$ must be isomorphic to some $L(\l)$. We now
prove a quantum Clebsch-Gordan formula for the finite dimensional
simple $U$-modules.

\begin{theorem}
Let $\l=(\l_1,\l_2), \mu=(\mu_1,\mu_2)\in (k^{\ti})^2$ with
$\l_1=\ep_1q^m,\, \mu_1=\ep_2q^n$ for some $m, n\>0$, where
$\ep_1=\pm 1$ and $\ep_2=\pm 1$. Then there exists an isomorphism of
$U$-modules
$$L(\l)\ot L(\mu)\simeq \bigoplus_{i=0}^{\min\{m,n\}}L(\eta^{(i)}),$$
where $\eta^{(i)}=(\ep_1\ep_2q^{m+n-2i},q^i\l_2\mu_2)$.
\end{theorem}
\begin{proof}
Similar to the proof of \cite[VII7.1]{kas}.
\end{proof}

Denote by $\mcf$ the category of all finite dimensional modules over
$U$. Note that the category $\mcf$ is not semisimple, see the
following example.
\begin{example}
For $\ep_1,\ep_2 \in \{1,-1\}$ and $a \in k^{\ti}$, let
$V(\ep_1,\ep_2,a)$ be a vector space with a basis $\{v_1,v_2\}$ and
define an action of $U$ as follows:
       \begin{eqnarray*}
&&K_1v_i=\ep_iv_i,\,K_2v_1=av_1,\,K_2v_2=aq^{-2}v_2,\\
       && E_1v_i=E_2v_i=0,\\
        &&F_1v_i=F_2v_2=0, \,F_2v_1=v_2.\\
 \end{eqnarray*}
Then $V(\ep_1,\ep_2,a)$ is a $2$-dimensional indecomposable module,
which is not semisimple.
\end{example}

The module structure of $\ml$ is complicated and is not well
understood. We shall make an attempt in this direction.

Let $u\in \ml$ with $E_1u=0$. Then $u$ is a linear combination of
the following vectors:
\begin{itemize}
  \item $F_2^{t_3}v_{t_2}$;
  \item $F_1^{m_1+t_3-t_2+1}F_2^{t_3}v_{t_2},$ if $\l(K_1)=\pm
  q^{m_1}$ for some $m_1 \in \Z$ and $m_1+t_3-t_2\>0.$
\end{itemize}

Suppose now that $\l(K_1)=\pm  q^{m_1}\text{ for some }m_1 \in \Z$.
Define $N(\l)$ be the subspace of $\ml$ spanned by the vectors
$F_1^{t_1}F_2^{t_3}v_{t_2}$, where $t_1\>m_1+t_3-t_2+1$ if
$t_2-t_3\<m_1$ and $t_1\> 0$ otherwise. Then
$N(\l)=\text{span}\{F_1^{t_1}F_2^{t_3}v_{t_2}|t_1>\max\{t_3-t_2+m_1,-1\}~\}$

Let $\Pi (M(\l))$ be the weight set of $M(\l)$. If $\mu \in \Pi
(M(\l))$, then $\mu =(q^{-2i+j}\l_1,q^{-2j+i}\l_2)$ for some $i,j\in
\N_0$. Note that $q^{s-2t}=q^{t-2s}=1$ for $s,t\>0$ if and only if
$s=t=0$. Hence $i,j$ are determined uniquely by $\mu$. We also have
that $\{F_1^iF_2^jv,F_1^{i-1}F_2^{j-1}v_1,\cds,
F_1^{i-s}F_2^{j-s}v_s\}$ is a basis of $\ml _{\mu}$ from the
decomposition (\ref{24}), where $s=\min\{i,j\}$.

\begin{remark}
 $N(\l)$ is not a submodule of $\ml$. For example, take $t_2,t_3\in \N_0$ with
 $m_1+t_3-t_2+1=0$. Then
$F_2^{t_3}v_{t_2}\in N(\l)$ but $F_2F_2^{t_3}v_{t_2}\notin  N(\l)$.
\end{remark}

\section{\bf The case $q$ is a root of unit}\selabel{3}
In this section, assume that $q$ is a primitive $l$-th root of
unity. Let $l'=l$ if $l$ is odd and $l'=l/2$ if $l$ is even. Then
$[l']_q=0$. Let $\zeta _l$ denote the set of all the $l-$th roots of
unity in $k$.
\subsection{The finite dimensional quotient Hopf algebra $\su$}
In the sequel, we shall construct a finite dimensional Hopf algebra
from $U$.
\begin{lemma}\lelabel{231}
The elements
$E_1^l,E_{1,2}^l,E_2^l,F_1^l,F_{1,2}^l,F_2^l,K_1^l,K_2^l$ are in the
center of $U$.
\end{lemma}
\begin{proof}
It follows from the relations $(10)-(15)$ and $(17)-(22)$.
\end{proof}

Let $I$ be the ideal of $U$ generated by the first six elements in
\leref{231} and the elements $K_i^l-1(i=1,2)$. Then $I$ is a Hopf
ideal. Define $\su\coloneqq U/I$. This is a finite dimensional Hopf
algebra. As an algebra, $ \su$ is generated by generators
$E_i,F_i,K_i(i=1,2)$ subject to the relations $(10)-(15)$ and the
following relations:
\begin{align}\eqlabel{23}E_i^l=0=F_i^l,\quad
K_i^l=1 \text{ for } i=1,2.\end{align}

Let $I^i=I\cap U^i$ for $i\in \{+,-,0\}$. Then $I^i$ is an ideal of
$ U^i$ and $I\simeq I^p\ot I^q\ot I^r$ as vector spaces, where
$(p,q,r)$ is a permutation of $(+,-,0)$. So $I=I^p  I^q I^r$. Let
$\su ^i=(U^i+I)/I$. By virtue of \eqref{16}, we have
\begin{proposition}\prlabel{232}
$$\su \simeq \su ^p\ot \su ^q \ot \su ^r.$$
\end{proposition}

By abuse of language, we denote the images of the elements $x\in U$
in $\su $ (resp. $\su^+,\su^-,\su ^0$) under the natural map again
by $x$. Then we have the following corollary about the structure of
$\su$.

\begin{corollary}\colabel{233}
$\su =\su ^p \su ^q  \su ^r$ and has a $k$-basis
$\{E_1^{s_1}E_{1,2}^{s_2}E_2^{s_3}K_1^{r_1}K_2^{r_2}F_1^{t_1}F_{1,2}^{t_2}F_2^{t_3}
|0\< s_i,t_i ,r_j \<l-1,1\<i\<3,1\<j\<2\}$.
\end{corollary}

Denote the Borel subalgebra $\su ^+\su ^0 (resp. \su ^0 \su ^-)$ by
$\su ^{\>0}(resp. \su ^{\<0})$. For any $0\< s_i,r_j
\<l-1,1\<i\<3,1\<j\<2$, let
$$deg_{\su^{\>0}}(E_1^{s_1}E_{1,2}^{s_2}E_2^{s_3}K_1^{r_1}K_2^{r_2})=
(s_1+s_2)\a_1+(s_2+s_3)\a_2=deg_{\su^{\<0}}(K_1^{r_1}K_2^{r_2}F_1^{s_1}F_{1,2}^{s_2}F_2^{s_3}).$$
Then $\su ^{\>0}$ and $\su ^{\<0}$ are both $Q$-graded Hopf algebras
with the gradings $deg_{\su^{\>0}}$ and $deg_{\su^{\<0}}$,
respectively.

Then there is a unique skew Hopf pairing $\vf:\su ^{\>0}\ot \su
^{\<0}\lra k$ such that
\begin{align*}
&\vf(1,1)=\vf (1,K_i)=\vf (K_i,1)=1,\\
&\vf (x,y)=0 \text{   if }x,y
 \text{ are monomials with } deg_{\su^{\>0}}(x)\neq deg_{\su^{\<0}}(y),\\
&\vf (E_i,F_j)=\d _{i,j}\d_ {i,1} \frac{1}{q^2-1},\\
&\vf (K_i,K_j)=q^{a_{ij}},\vf (K_i,K_j^{-1})=q^{-a_{ij}},
\end{align*}
for $1\<i,j\<2$. Thus one can turn $D_{\vf}(\su^{\>0},\su
^{\<0})=\su^{\>0}\ot\su ^{\<0}$ into a Hopf $k$-algebra as in
\seref{1}.

For convenience, we write $D_{\vf}$ instead of
$D_{\vf}(\su^{\>0},\su ^{\<0})$. Then $D_{\vf}$ can be described as
follows.

$D_{\vf}$ is generated, as an algebra, by $E_i,F_i,K_i, K_i^{-1},
\widetilde{K}_i,\widetilde{K}_i^{-1}(1\<i\<2)$ subject to the
relations:
\begin{align*}
&K_iK_j=K_jK_i,\quad K_iK_i^{-1}=K_i^{-1}K_i=1, \\
&\widetilde{K}_i\widetilde{K}_j=\widetilde{K}_j\widetilde{K}_i,\quad
\widetilde{K}_i\widetilde{K}_i^{-1}=\widetilde{K}_i^{-1}
\widetilde{K}_i=1,\quad K_i\widetilde{K}_j=\widetilde{K}_jK_i ,\\
&K_iE_jK_i^{-1}=q^{a_{ji}}E_j,\quad
K_iF_jK_i^{-1}=q^{-a_{ji}}F_j,\\
&\widetilde{K}_iE_j\widetilde{K}_i^{-1}=q^{a_{ji}}E_j,\quad
\widetilde{K}_iF_j\widetilde{K}_i^{-1}=q^{-a_{ji}}F_j,\\
&E_iF_j-F_jE_i=\d_{ij}\d_{i1}(K_i-\widetilde{K}_i^{-1})/(q-q^{-1})\\
&E_i^l=0=F_i^l,\quad K_i^l=1=\widetilde{K}_i^l.
\end{align*}
The coalgebra structure is given by
\begin{eqnarray*}
&\D(E_i)=K_i\ot E_i+E_i\ot 1,&\D(F_i)=1\ot F_i+F_i\ot \widetilde{K}_i^{-1},\\
&\D(K_i)=K_i\ot K_i,&\D(\widetilde{K}_i)=\widetilde{K}_i\ot \widetilde{K}_i,\\
&S(E_i)=-K_i^{-1}E_i,&S(F_i)=-F_i\widetilde{K}_i,\quad S(K_i)=K_i^{-1},
S(\widetilde{K}_i)=\widetilde{K}_i^{-1},\\
& \ep (E_i)=0=\ep (F_i), &\ep (K_i)=1=\ep (\widetilde{K}_i) .
\end{eqnarray*}
This is similar to the Drinfeld quantum double of $\su^{\<0}$. In
fact, $D_{\vf}$ is quasi-isomorphic to the Drinfeld quantum double
$D(\su _q^{\>0}(\mathfrak{sl}(3)),\su _q^{\<0}(\mathfrak{sl}(3)))$,
i.e., the categories of their comodules are monoidally equivalent
(\cite{did2}).

Recall that Andruskiewitsch and Schneider gave a classification of
finite dimensional pointed Hopf algebras, see \thref{1.20}. Since
$D_{\vf}$ is a finite dimensional pointed Hopf algebra, one can
obtained it by the general construction method (Similarly to the
process to construct $\su$ before.): choose a datum and a linking
datum; define a big Hopf algebra; then modulo some suitable Hopf
ideal such that the induced quotient is isomorphic to $D_{\vf}$.
Here the two data of $D_{\vf}$ are the same with those of $\su$
except abelian group $\G$. Thus we have
\begin{lemma}\lelabel{234}
The map which sends $E_i$ to $E_i$, $F_i$ to $F_i$, $K_i^{\pm 1}$ to
$K_i^{\pm 1}$ and $\widetilde{K}_i^{\pm 1}$ to $K_i^{\pm 1}$ can be
extended uniquely to a surjective Hopf algebra homomorphism $\pi
:D_{\vf} \lra \su$. Moreover, $\ker \pi =\langle
K_1-\widetilde{K}_1, K_2-\widetilde{K}_2\rangle$.
\end{lemma}
\begin{proof}
It is easy to check.
\end{proof}

The categories of modules over $D_{\vf}$ and $\su$ are connected
closely.

Similar to the definition of the weight modules over $U$, a
$D_{\vf}$(resp. $\su$)-module is called a weight module if it can be
decomposed as a direct sum of $1$-dimensional simple modules over
$D_{\vf}^0 $(resp. $\su$), where $D_{\vf}^0 $ is the subalgebra of
$D_{\vf} $ generated by $K_i, \widetilde{K}_i,(1\<i\<2)$.
\begin{lemma}\lelabel{235}
Every $D_{\vf}$(resp. $\su$)-module is a weight module.
\end{lemma}
\begin{proof}
It follows from the fact that the subalgebra generated by $
K_i,\widetilde{K_i}$(resp. $ K_i $), $i=1,2$, is a group algebra of
finite abelian group over an algebraically closed field $k$ of
char$(k)=0$.
\end{proof}

From the relations defining $D_{\vf}$, we know that $K_i\wt
{K}_i^{-1}$ are in the center of $D_{\vf}$. Furthermore, we have
\begin{lemma}\lelabel{236}
Let $M$ be an indecomposable $D_{\vf}$-module. Then $K_i\wt
{K}_i^{-1}$ acts on $M$ as some scalar $z_i\in \zeta _l$, $i=1,2$.
\end{lemma}
\begin{proof}
By \leref{235}, $K_i,\wt {K}_i^{-1}$ act semisimply on $M$. Since
$K_i$ commutes with $\wt {K}_i^{-1}$, the action of $K_i\wt
{K}_i^{-1}$ is semisimple on $M$, i.e., there is a direct sum
decomposition $M=\op _{z \in k^2}M_z$, where $M_z=\{m\in M|K_i\wt
{K}_i^{-1}m=z_im,i=1,2\}$. Indeed, $z=(z_1,z_2) \in \zeta _{l}\ti
\zeta _l$ since $(K_i\wt {K}_i^{-1})^l=1$. Then the claim follows
from that $K_i\wt {K}_i^{-1}$ is in the center of $D_{\vf} $ and $M$
is indecomposable.
\end{proof}

In the rest of this subsection, we assume that $l$ is odd. For any
$z_1\in \zeta_{l}$, fix an element $z_1^{1/2}\in k$ such that
$(z_1^{1/2})^2=z_1$. In particular, let $1^{1/2}=1$.

Let $M$ be a $D_{\vf}$-module and $z=(z_1,z_2)\in \zeta _{l}\ti
\zeta _l$ such that $K_i\wt {K}_i^{-1}$ acts on $M$ as the scalar
$z_i$. Then there is a $k$-algebra homomorphism $\pi _z: D_{\vf}\lra
\su$ such that
\begin{align*}
&\pi_z (E_1)=z_1^{1/2}E_1,\,\pi _z(F_1)=F_1,\,\pi_z
(K_1)=z_1^{1/2}K_1,\,\pi_z (\wt{K}_1)=z_1^{-1/2}K_1,\\
&\pi_z (E_2)=E_2,\,\pi _z(F_2)=F_2,\,\pi_z (K_2)=K_2,\,\pi_z
(\wt{K}_2)=z_2^{-1}K_2.
\end{align*}

One can easily check that $\pi _z$ is well defined and the kernel of
$\pi _z$, which is the ideal generated by
$K_i\wt{K}_i^{-1}-z_i(i=1,2)$, annihilates the module $M$. Thus $M$
becomes a $\su$-module through the homomorphism $\pi _z$.
\begin{lemma}\lelabel{237}
Every indecomposable $D_{\vf}$-module is the pull-back of some
$\su$-module through an algebra homomorphism $\pi _z$.
\end{lemma}

Let $M$ be a $\su$-module and $z=(z_1,z_2)\in \zeta_{l}\ti \zeta_l$.
Denote by $M_z$ the pull-back of $M$ through the algebra
homomorphism $\pi _z$. Let $\ep_z$ be the $1$-dimensional
representation of $D_{\vf}$ defined by
$$\ep_z(E_i)=0=\ep_z(F_i),\, \ep_z(K_1)=z_1^{1/2},\,
\ep_z(\wt{K}_1)=z_1^{-1/2},\,\ep_z(K_2)=1,\,\ep_z(\wt{K}_2)=z_2^{-1}.$$

One can easily check that $\ep_z$ is well-defined.

\begin{lemma}\lelabel{238}
Let $z\in \zeta_{l}\ti \zeta_l$ and $M$ be a module over $\su$. Then
$$M_z\simeq \ep_z\ot M_1,\text{ where }
1=(1,1).$$
\end{lemma}
\begin{proof}
It follows from a direct verification.
\end{proof}

From \leref{238}, we obtain the following theorem which illustrates
the relation between the module categories of $D_{\vf}$ and $\su$.
\begin{theorem}\thlabel{239}
The category $_{D_{\vf}} \mcm$ is equivalent to the direct product
of $|\zeta_{l}\ti \zeta_l|$ copies of the category $_{\su} \mcm$.
\end{theorem}
\begin{remark}
\begin{enumerate}
  \item For the case $l$ is even, one can also define an algebra
homomorphism $\pi_z$ for any given $z=(z_1,z_2)\in \zeta_l \ti
\zeta_l$ if $z_1^{1/2}\in \zeta_l$.
  \item One can define a skew Hopf pairing $\psi:U^{\>0}\ot U^{\<0}\lra
  k$ and form the corresponding Hopf algebra $D_{\psi}(U^{\>0}, U^{\<0})$, then one can discuss the
  relation between the categories of weight modules over $D_{\psi}(U^{\>0},
  U^{\<0})$ and $U$, which is similar to the discussion in \cite{huz} for
  the quantized enveloping algebra.
\end{enumerate}
\end{remark}

\subsection{The representation theory over $\su$}
 Now let us concentrate on the representation theory of $\su$.
First we list all the simple modules over $\su$.

For any $0\<m_1,m_2<l$, define an $m_1+1$ dimensional $\su$-module
$V=span\{w_0,w_1,\cds, w_{m_1}\}$ as follows:
\begin{align*}
&K_1w_j=q^{m_1-2j}w_j,\quad &&K_2w_j=q^{m_2+j}w_j,\\
&E_1w_j=[m_1+1-j]w_{j-1},\quad &&F_1w_j=[j+1]w_{j+1},\\
&E_2w_j=F_2w_j=0,
\end{align*}
where $0\<j\<m_1$ and $w_{-1}=0=w_{m_1+1}$. Then one can easily
check that $V$ becomes a $\su$-module, denoted by $V(m_1,m_2)$, and
that $V(m_1,m_2)$ is simple. Furthermore, for any $0\<n_1,n_2<l$,
$V(m_1,m_2)\simeq V(n_1,n_2)$ if and only if $m_1=n_1,m_2=n_2$.
\begin{proposition}\prlabel{2310}
Every simple $\su$-module is isomorphic to some $V(m_1,m_2)$.
\end{proposition}
\begin{proof}
It can be shown by a direct verification. It also follows from
\prref{2112} and the fact that $\su$ is a quotient of $U$.
\end{proof}

Let $P(m_1,m_2)$ be the projective cover of $V(m_1,m_2)$.

For any pair $(m_1,m_2)$ of integers with $0\<m_1,m_2\<l-1$, let
$V=kv$ be a $1$-dimensional $\su ^{\>0}$-module with the action
given by $E_1\cd v =0=E_2\cd v$, $K_1\cd v=q^{m_1}v$ and $K_2\cd
v=q^{m_2}v$. Define the Verma module $M(m_1,m_2)=\su \ot
_{\su^{\>0}}V$.  Since $\su$ is a quotient of $U$, each $\su$-module
naturally becomes a $U$-module. As a $U$-module, $M(m_1,m_2)$ is
isomorphic to $M(\l)/(I\cd w)$, where $\l=(q^{m_1},q^{m_2})$ and
$w\neq 0$ is a weight vector in $M(\l)$ with weight $\l$ such that
$M(\l)=U\cd w$. By \prref{212}, $M(m_1,m_2)$ has a unique maximal
submodule as $U$-module. Consequently, $M(m_1,m_2)$ has a unique
maximal submodule as $\su$-module.

\begin{proposition}\prlabel{2311}
For any pair $(m_1,m_2)$ of integers with $0\<m_1,m_2\<l-1$,
$M(m_1,m_2)$ is isomorphic to some quotient of $P(m_1,m_2)$.
\end{proposition}
\begin{proof}
Let $f:M(m_1,m_2)\ra V(m_1,m_2)$, $g:P(m_1,m_2)\ra V(m_1,m_2)$ be
the canonical projections. By the projectivity of $P(m_1,m_2)$,
there is a homomorphism $g':P(m_1,m_2)\ra M(m_1,m_2)$ such that
$fg'=g$. Take a weight vector $v\in P(m_1,m_2)$ such that $g(v)\neq
0$. Without loss of generality, one may assume that $g(v)$ is a
highest weight vector by a suitable choice of $v$. Since $\su
v\nsubseteq \ker g=\text{rad} P(m_1,m_2)$ and $\text{head}
P(m_1,m_2)\simeq V(m_1,m_2)$ is simple, $\su v+\text{rad}
P(m_1,m_2)=P(m_1,m_2)$. Then $\su v=P(m_1,m_2)$ by the Nakayama
Lemma. Similarly, one can show that $g'(v)\in M(m_1,m_2)$ is a
weight vector and $\su g'(v)=M(m_1,m_2)$. Thus $g'$ is surjective.
\end{proof}

Let $(m_1',m_2')$ be another pair of integers with
$0\<m_1',m_2'\<l-1$. Then it is easy to check that $\Hom
(P(m_1,m_2),P(m_1',m_2'))\neq 0$ if and only if $V(m_1,m_2)$ is a
composition factor of $P(m_1',m_2')$.

\begin{theorem}\thlabel{2312}
Assume that $(l,3)=1$. Then $\su$ is an indecomposable algebra.
\end{theorem}
\begin{proof}
By the discussion above, it suffices to show that each $V(m_1,m_2)$
is a composition factor of $P(0,0)$. If $V(m_1,m_2)$ is a
composition factor of $M(0,0)$, then it is a composition factor of
$P(0,0)$. However, $V(m_1,m_2)$ is a composition factor of $M(0,0)$
if and only if there exist $0\<t_2,t_3\<l-1$ such that
$\begin{pmatrix}
   -1 & 1 \\
   -1 & -2 \\
 \end{pmatrix}\begin{pmatrix}
                t_2 \\
                t_3 \\
              \end{pmatrix}\equiv \begin{pmatrix}
                                    m_1 \\
                                    m_2 \\
                                  \end{pmatrix}\quad (\text{ mod } l)$ since
                                  $F_2^{t_3}v_{t_2}(0\<t_2,t_3\<l-1)$
 (or their linear combination) are the only
 highest weight vectors.

 We claim
                                  that the equations
$$\begin{pmatrix}
   -1 & 1 \\
   -1 & -2 \\
 \end{pmatrix}\begin{pmatrix}
                t_2 \\
                t_3 \\
              \end{pmatrix}\equiv \begin{pmatrix}
                                    1 \\
                                    0 \\
                                  \end{pmatrix}\quad (\text{ mod } l) \text{ and }
\begin{pmatrix}
   -1 & 1 \\
   -1 & -2 \\
 \end{pmatrix}\begin{pmatrix}
                t_2 \\
                t_3 \\
              \end{pmatrix}\equiv\begin{pmatrix}
                                    0 \\
                                    1 \\
                                  \end{pmatrix}\quad (\text{ mod } l)$$have
                                  solutions in $\Z$.

In fact, the first equation is equivalent to
\begin{equation*}
\left\{ \begin{aligned}
        1 &= (2s+t)l+3(t_2+1) \\
                1  &=(t-s)l+3t_3
                          \end{aligned} \right.
                          \end{equation*}
for some $s,t\in \Z$. Since $(l,3)=1$, there exist $p,q\in \Z$ such
that $pl+3q=1$. Let $t_3=q,t_2=q-1,t=p,s=0$. Then this is a solution
of the equation set above. The second equation can be converted into
the equation set
\begin{equation*} \label{eq:2}
 \left\{ \begin{aligned}
        1 &= (-2s-t)l-3t_2 \\
                1  &=(s-t)l-3t_3
                          \end{aligned} \right.
                          \end{equation*}
for some $s,t\in \Z$. Let $t_2=-q=t_3, s=0, t=-p$. This is a
solution of the equation set above. So there exist
$t_2',t_3',t_2'',t_3''\in \Z$ such that
$$\begin{pmatrix}
   -1 & 1 \\
   -1 & -2 \\
 \end{pmatrix}\begin{pmatrix}
                t_2' \\
                t_3' \\
              \end{pmatrix}\equiv \begin{pmatrix}
                                    1 \\
                                    0 \\
                                  \end{pmatrix} \quad (\text{ mod }l) \text{ and }
\begin{pmatrix}
   -1 & 1 \\
   -1 & -2 \\
 \end{pmatrix}\begin{pmatrix}
                t_2'' \\
                t_3'' \\
              \end{pmatrix}\equiv \begin{pmatrix}
                                    0 \\
                                    1 \\
                                  \end{pmatrix} \quad (\text{ mod } l)$$
Let $t_2,t_3\in \Z$ with $0\<t_2,t_3\<l-1$ such that $t_2\equiv
m_1t_2'+m_2t_2''\quad (\text{ mod } l)$ and $t_3\equiv
m_1t_3'+m_2t_3''\quad (\text{ mod } l)$. Then
$$\begin{pmatrix}
   -1 & 1 \\
   -1 & -2 \\
 \end{pmatrix}\begin{pmatrix}
                t_2 \\
                t_3 \\
              \end{pmatrix}\equiv \begin{pmatrix}
                                    m_1 \\
                                    m_2 \\
                                  \end{pmatrix} \quad (\text{ mod } l)$$
\end{proof}

For the projective modules $P(m_1,m_2)$, we hope to give a more
detailed description. It is well known that there is a close
connection between the direct sum decomposition of the regular
module and the idempotent element decomposition of unit. Let us go
along this way.

For $1\<i,j\<l-1$, define $e_{i,j}=\frac{1}{l^2}\sum
_{s,t=0}^{l-1}q^{is+jt}K_1^sK_2^t$ in $\su$. Then we have
$e_{i,j}^2=e_{i,j},\, e_{i,j}e_{i',j'}=\d_{i,i'}\d_{j,j'}$ and
$\sum_{i,j=0}^{l-1}e_{i,j}=1$. So $\su =\op _{i,j=0}^{l-1}\su
e_{i,j}$ is a direct sum decomposition of $\su$ as regular module.
Each summand $\su e_{i,j}$ is a projective module.

\begin{proposition}\prlabel{2313}
For $1\<i,j\<l-1$, $\su e_{i,j}$ are non-isomorphic each other.
\end{proposition}
\begin{proof}
Assume that $\su e_{i,j}$ is isomorphic to $ \su e_{i',j'}$ for some
$(i',j')\neq (i,j)$, and let $f:\su e_{i,j}\rightarrow  \su
e_{i',j'}$ be an isomorphism. Then there are $a,b\in \su$ such that
$f(e_{i,j})=ae_{i',j'}$ and $f^{-1}(e_{i',j'})=be_{i,j}$, and so
$e_{i,j}=abe_{i,j}$. Note that $K_1e_{i,j}=q^{-i}e_{i,j},\,
K_2e_{i,j}=q^{-j}e_{ij}$. Hence $a,b\notin \su ^0$, and we may
assume that $a=\sum _{\textbf{s,t}\in
\Z_l^3}\a_{\textbf{s,t}}E^{\textbf{s}}F^{\textbf{t}},b=\sum
_{\textbf{s}',\textbf{t}'\in
\Z_l^3}\b_{\textbf{s}',\textbf{t}'}E^{\textbf{s}'}F^{\textbf{t}'}$,
where $\a_{\textbf{s,t}},\b_{\textbf{s}',\textbf{t}'} \in k$. It is
easy to see that $\{E^{\textbf{s}}F^{\textbf{t}}e_{ij}|
\textbf{s,t}\in \Z_l^3, i,j\in \Z_l\}$ is a $k$-basis of $\su$.
Hence we deduce $ab=1$ from $e_{i,j}=abe_{i,j}$. Comparing the
degree of terms on the both sides of $ab=1$, we have that $a=\a
_{0,0},b =\b _{0,0}$ in $k$. It contradicts with $a,b\notin \su ^0$.
This proved the claim.
\end{proof}

Recall that $P(m_1,m_2)$ is the projective cover of $V(m_1,m_2)$ and
$\dim V(m_1,m_2)=m_1+1$. From the representation theory of finite
dimensional algebra, we know that
\begin{equation}\label{25}\su \simeq \bigoplus
_{0\<m_1,m_2\<l-1}(m_1+1)P(m_1,m_2).\end{equation}

Thus $\su e_{ij}$ may not be indecomposable, i.e., $e_{ij}$ may not
be primitive. It is difficult to give all the primitive idempotent
elements in $\su$. In the following we will turn this into the
question to solve some equation set.

Now we assume that $l$ is odd in the sequel.

Let $\su _1$ be the subalgebra of $\su$ generated by $E_1,F_1,K_1$.
Then $\su _1$ is isomorphic to the quotient of quantum group
$U_q(\mathfrak{sl}(2))$ modulo the ideal $\langle E^l,F^l,
K^l-1\rangle$. Thus $\su_1$ just is the so-called ``small'' quantum
group(cf. \cite{xia}).

From the results in \cite{xia}, we know that $\su_1\simeq \op
_{i=0}^{l-1}(i+1)P_i$ as $\su _1$-modules, where each $P_i$ is
indecomposable and non-isomorphic. Thus there is a family of
primitive orthogonal idempotent elements $f_1,f_2,\cds,
f_N(N=l(l+1)/2)$ in $\su_1$ such that $1=\sum f_i$.

Define $e_j(K_2)=\frac{1}{l}\sum _{i=0}^{l-1}(q^jK_2)^i$. Then
$1=\sum _{j=0}^{l-1}e_j(K_2)$ is a decomposition of orthogonal
idempotent elements. Obviously, $f_ie_j(K_2)$ is an idempotent
element, denoted by $f_{ij}$. One can check that $f_{ij}$ are
pairwise orthogonal and $1=\sum_{i,j} f_{ij}$. By virtue of the
decomposition (\ref{25}) and the Krull-Schmidt theorem, $f_{ij}$ is
a primitive idempotent element in $\su$.

Thus we only need to give all primitive idempotent elements in
$\su_1$ in order to give all the primitive idempotent elements in
$\su$.

For convenience, write $K,E,F,e_i$ for $K_1,E_1,F_1$ and
$e_i(K_1)=\frac{1}{l}\sum _{s=0}^{l-1}(q^iK_1)^s$ respectively.
Clearly, $\{e_iE^sF^t|0\<i,s,t\<l-1\}$ is a basis of $\su_1$.

Assume that $e_i$ is not primitive. Then there is an idempotent
element decomposition $e_i=\sum_r e_{i,r}$, where $e_{i,r}$ are
orthogonal primitive elements. Since $\su =\op _{i=0}^{l-1}e_i\su
=\op _{i=0}^{l-1}e_i\su ^+\su^-$, we have $e_{i,r}\in
e_i\su^+\su^-,$ and so $e_{i,r}=\sum _{t=1-l}^{l-1}e_iP_{r,t}$,
where $P_{r,t}=\sum _{s=0}^{l-1}a_{r,t,s}E^sF^{t+s}$ and $F^{t+s}=0$
for $t+s<0$ or $t+s>l-1$. Note that $e_i=\sum _{r,t}e_iP_{r,t}=\sum
_{r,t}P_{r,t}e_{i-2t}$. Multiply $e_i$ on the equation from the both
sides one gets $e_i=\sum_rP_{r,0}e_i=\sum_re_iP_{r,0}$, here we use
the assumption that $l$ is odd. Hence $P_{r,t}=0$ for $t\neq 0$.
This means that $e_{i,r}=e_iP_{r,0}$, where
$P_{r,0}=\sum_{s=0}^{l-1}a_{r,s}E^sF^s$.

From the discussion above, we only need to give all the idempotent
elements of form $\sum_{s=0}^{l-1}a_se_iE^sF^s$.

\begin{lemma}\lelabel{2314}
$$F^mE^s=\sum _{j=0}^{\min\{m,s\}}[j]_q!\mat{m}{j}_q\mat{s}{j}_qE^{s-j}
\prod _{r=j-m-s+1}^{2j-m-s}[K^{-1};r]_qF^{m-j}$$
\end{lemma}
\begin{proof}
Define an automorphism $\vf$ of $U_{\upsilon}(\mathfrak{sl}(2))$ by
$\vf E=F,\vf F=E,\vf K=K,\vf \upsilon=\upsilon^{-1},$ where
$\upsilon$ is an indeterminate over $k$. Then it follows from
applying $\vf$ to the identity in \cite[Lemma1.7]{kas} and
$U_q(\mathfrak{sl}(2))=U_{\upsilon}(\mathfrak{sl}(2))/(\upsilon-q)$.
\end{proof}

Then we have the following identity:
$$e_iE^mF^me_iE^sF^s=\sum _{j=0}^{\min\{m,s\}}[j]_q!^2\mat{m}{j}_q\mat{s}{j}_q
\mat{m+s+i}{j}_qe_iE^{m+s-j} F^{m+s-j}.$$

Let $e_{i,r}=\sum_{p=0}^{l-1}a_pe_iE^pF^p$ for a fix $r$. By using
the identity above repeatedly, we have
\begin{theorem}\thlabel{2315}
For some fixed $i\in \Z_l$, $e_{i,r}$ is an idempotent element if
and only if \\$(a_0,a_1,\cds, a_{l-1})$ is a solution of the
equation set$$a_p=\sum
_{\substack{0\<m,s\<l-1,\\0\<j\<\min\{m,s\},\\m+s-j=p}}a_{m,s,j}a_ma_s
,\quad 0\<p\<l-1,$$ where
$a_{m,s,j}=[j]_q!^2\mat{m}{j}_q\mat{s}{j}_q \mat{m+s+i}{j}_q$.
\end{theorem}

Let $[n]$ be the residual class of $n$ modulo $l$ and $\ov{n}$ be
the minimal non-negative integer in $[n]$. It is easy to see that:
$a_{m,s,j}=a_{s,m,j}$ and that $\mat{m+s+i}{j}_q=0$ if and only if
$\ov{m+s+i}<j$. If $(a_0,a_1,\cds, a_{l-1})$ is a solution of the
equation set, then $a_0=0$ or $a_0=1$. Moreover, $(1-a_0,-a_1,\cds,
-a_{l-1})$ is also a solution of the equation set.
\begin{remark}\relabel{2316}
The equation set in \thref{2315} is solvable recursively. For the
special cases $l=3$, we will give all the solutions in the below.
This result is new for the small quantum group $\su_1$.
\end{remark}

\begin{example}
When $l=3$, we can get a decomposition of orthogonal primitive
idempotents of unity through solving  the equation set in
\thref{2315}.
\begin{equation*}
\begin{split}
1=&(-e_0EF+e_0E^2F^2)+(e_0+e_0EF-e_0E^2F^2)\\&+(e_1-e_1E^2F^2)
+e_1E^2F^2\\&+(e_2EF-e_2E^2F^2)+(e_2-e_2EF+e_2E^2F^2),\end{split}
\end{equation*}
where $e_i=\frac{1}{3}\sum_{j=0}^2(q^iK)^j\in \su_1$ for $i=0,1,2$.
\end{example}

\end{document}